\documentclass{article}
\usepackage{spconf,amsmath,graphicx}
\usepackage{relsize}
\usepackage[T1]{fontenc}
\usepackage{microtype}
\usepackage{lscape}
\usepackage{etoolbox}
\usepackage{adjustbox}
\usepackage{graphicx}
\usepackage{booktabs} 
\usepackage{comment}
\usepackage{csquotes}
\usepackage{enumerate}
\usepackage{algorithm}
\usepackage{algorithmic}
\usepackage{comment}

\usepackage{thmtools}
\usepackage{thm-restate}
\usepackage{caption}
\usepackage{subcaption}
\usepackage{url}

\usepackage{amsmath}
\usepackage{amssymb}
\usepackage{mathtools}
\usepackage{amsthm}
\usepackage{siunitx}
\usepackage[backref=page]{hyperref}
\usepackage{cleveref}

\newtheorem{theorem}{Theorem}[section]
\newtheorem*{theorem*}{Theorem}

\theoremstyle{definition}
\newtheorem{definition}[theorem]{Definition}

\theoremstyle{remark}

\usepackage[textsize=tiny]{todonotes}

\newcommand{\subheading}[1]{\smallskip\noindent\textbf{#1.} }

\newtoggle{arxiv}
\togglefalse{arxiv}
\newcommand{\tSC}{\textsmaller{SC}}
\newcommand{\UpL}[1]{L_{#1}^\text{up}}
\newcommand{\DownL}[1]{L_{#1}^\text{down}}

\newcommand{\Z}{\mathbb{Z}}

\newcommand{\SC}{\mathcal{S}}

\newcommand{\R}{\mathbb{R}}

\DeclareMathOperator{\Ima}{Im}

\everymath{\medmuskip=2.5mu minus 2.5mu\thickmuskip=2.5mu minus 2.5mu}

\newcommand\michael[1]{\noindent{\textcolor{magenta}{[MTS: #1]}}}
\newcommand\vincent[1]{\noindent{\textcolor{magenta}{[VPG: #1]}}}
\renewcommand\vincent[1]{}\renewcommand\michael[1]{}

\title{Disentangling the spectral properties of the Hodge Laplacian:\\ Not all small eigenvalues are equal}

\name{Vincent P. Grande\thanks{\textsmaller{VPG} acknowledges funding by the German Research Council (\textsmaller{DFG}) within Research Training Group 2236 (\textsmaller{U}n\textsmaller{RAV}e\textsmaller{L}).}, Michael T. Schaub\thanks{\textsmaller{MTS} acknowledges partial funding by the Ministry of Culture and Science (\textsmaller{MKW}) of the German State of North Rhine-Westphalia ("\textsmaller{NRW} R\"uckkehrprogramm") and the European Union (\textsmaller{ERC}, \textsmaller{HIGH-HOPeS}, 101039827). Views and opinions expressed are however those of the author(s) only and do not necessarily reflect those of the European Union or the European Research Council Executive Agency. Neither the European Union nor the granting authority can be held responsible for them.}}
\address{Department for Computer Science,
	\textsmaller{RWTH} Aachen University,
	Germany
}
\begin{document}
\ninept
\maketitle
\begin{abstract}
    The rich spectral information of the graph Laplacian has been instrumental in graph theory, machine learning, and graph signal processing for applications such as graph classification, clustering, or eigenmode analysis.
	Recently, the Hodge Laplacian has come into focus as a generalisation of the ordinary Laplacian for higher-order graph models such as simplicial and cellular complexes.
	Akin to the traditional analysis of graph Laplacians, many authors analyse the smallest eigenvalues of the Hodge Laplacian, which are connected to important topological properties such as homology.
	However, small eigenvalues of the Hodge Laplacian can carry different information depending on whether they are related to curl or gradient eigenmodes, and thus may not be comparable.
	We therefore introduce the notion of persistent eigenvector similarity and provide a method to track individual harmonic, curl, and gradient eigenvectors/-values through the so-called persistence filtration, leveraging the full information contained in the Hodge-Laplacian spectrum across all possible scales of a point cloud.
	Finally, we use our insights (a) to introduce a novel form of Hodge spectral clustering
	and (b) to classify edges and higher-order simplices based on their relationship to the smallest harmonic, curl, and gradient eigenvectors.
\end{abstract}
\begin{keywords}
Spectral Clustering, Topological Signal Processing, Persistent Homology, Topological Data Analysis
\end{keywords}
\section{Introduction}
In many areas of Mathematics, Physics, and Signal Processing, the eigenvalues of a linear operator or matrix encode important properties of the system governed by the operator.
Accordingly, harmonic and spectral analysis approaches have been proven to be useful in many application areas.
In the context of the analysis of networks and graphs, the idea of exploiting the spectral properties of an algebraic representation of the system is deeply engrained in spectral graph theory~\cite{godsil2001algebraic, chung1997spectral}.
While the initial focus of spectral graph theory has historically been on the adjacency matrix, the combinatorial and normalized Laplacian matrices (and variants thereof) have replaced the adjacency matrix as the operator typically considered, which may be explained by their favourable structure (e.g., positive semidefiniteness) and their deep connections to continuous mathematics.
Intuitively, the Graph Laplacian is an operator that encodes a notion of \enquote{local difference} between a node and its neighbours and may be seen as a discrete equivalent of the continous Laplacian operator.
Its $0$-eigenvalues encode connected components, whereas remaining smaller eigenvalues carrying information on graph connectivity \cite{spielman2012spectral} are used for spectral clustering \cite{vonLuxburg:2007} and graph classification and matching~\cite{fiori2015spectral}.
Further examples of spectral characteristics that are informative for applications include the spectral radius of the adjacency matrix, which governs the epidemic spreading threshold~\cite{saha2015approximation}, the second eigenvalue of the Laplacian, which characterizes the convergence of consensus dynamics~\cite{aysal2008accelerated}, and the spectral gap of the transition matrix, which is paramount for the  mixing behaviour of Markov chains~\cite{boyd2004fastest}.

In these archetypal examples, the relevant feature of an eigenvalue is its numerical value; eigenvalues are ``equal'' in other regards.
In this paper we will argue that this mantra does not hold for the spectrum of the Hodge Laplacian, a higher order generalisation of the Graph Laplacian for simplicial complexes (\textsmaller{SC}s).
\textsmaller{SC}s are higher-order generalisations of graphs which allow to incorporate complex relationships and geometric and topological information beyond pairwise interactions.
\textsmaller{SC}s have enjoyed a renaissance in the context of topological signal processing in the last years \cite{Barbarossa:2020, Schaub:2021}.
Importantly, the set of eigenvector/eigenvalue pairs of the Hodge Laplacian consists of three distinct types of eigenvectors: \emph{harmonic}, \emph{gradient}, and \emph{curl} eigenvectors, which correspond to the Hodge decomposition of the simplicial signal spaces.
Accordingly, rather than treating the spectrum of the Hodge Laplacian simply as a set of numbers whose numerical values are of interest, in this paper we more explicitly consider the gradient, curl and harmonic eigenvalues as separate to leverage the rich additional information provided by the Hodge Laplacian spectrum. 

\subheading{Contributions and Outline}
In \Cref{sec:background}, we give a brief overview over the concepts of topological signal processing, simplicial complexes, $\alpha$-filtrations, and persistent Laplacians.
In \Cref{sec:manyeigenvectors}, we track and analyse the evolution of the eigenvectors of the persistent Laplacian across the stages of an $\alpha$-filtrations in theory and experiments.
We introduce the concepts of persistent eigenvector similarity and matchings to give a new form of persistent eigenvalue diagrams.
We apply these insights on data in \Cref{sec:applications}, developing the notions of \emph{Hodge spectral clustering} (\Cref{fig:SpectralClustering}) and the \textsmaller{HGC}-values for inferring simplex 
roles (\Cref{fig:edgeroles}).
The accompanying code can be found \href{https://git.rwth-aachen.de/netsci/2024-disentangling-the-spectral-properties-of-the-hodge-laplacian}{here}.

\subheading{Related Work}
The persistent Laplacians were first defined in \cite{wang2020persistent}.
However, the authors did not distinguish between the gradient and curl eigenvalues in their analysis.
Further work on the persistent Laplacians includes \cite{memoli2022persistent, Davies2023, liu2023algebraic}.
The Hodge Laplacian has already been 
studied in the topological signal processing literature \cite{Schaub:2021}.
Ideas of spectral clustering \cite{vonLuxburg:2007} have been applied in some form to simplicial complexes.
However, Refs.\ \cite{Ebli2019, Grande:2023} only consider harmonic eigenvectors, Ref.\ \cite{Krishnagopal2021} is only interested in the support of the eigenvectors and disregards the gradient/curl differences, and Refs.\ \cite{reddy2023clustering,govek2019clustering} are only interested in point clustering.
The gradient and curl spaces of the Hodge Laplacian have been explored in \cite{schaub2014structure}
without the analysis of the connection to filtrations of simplicial complexes.
	\begin{figure*}[tb!]
		\begin{center}
			\includegraphics[width=\linewidth]{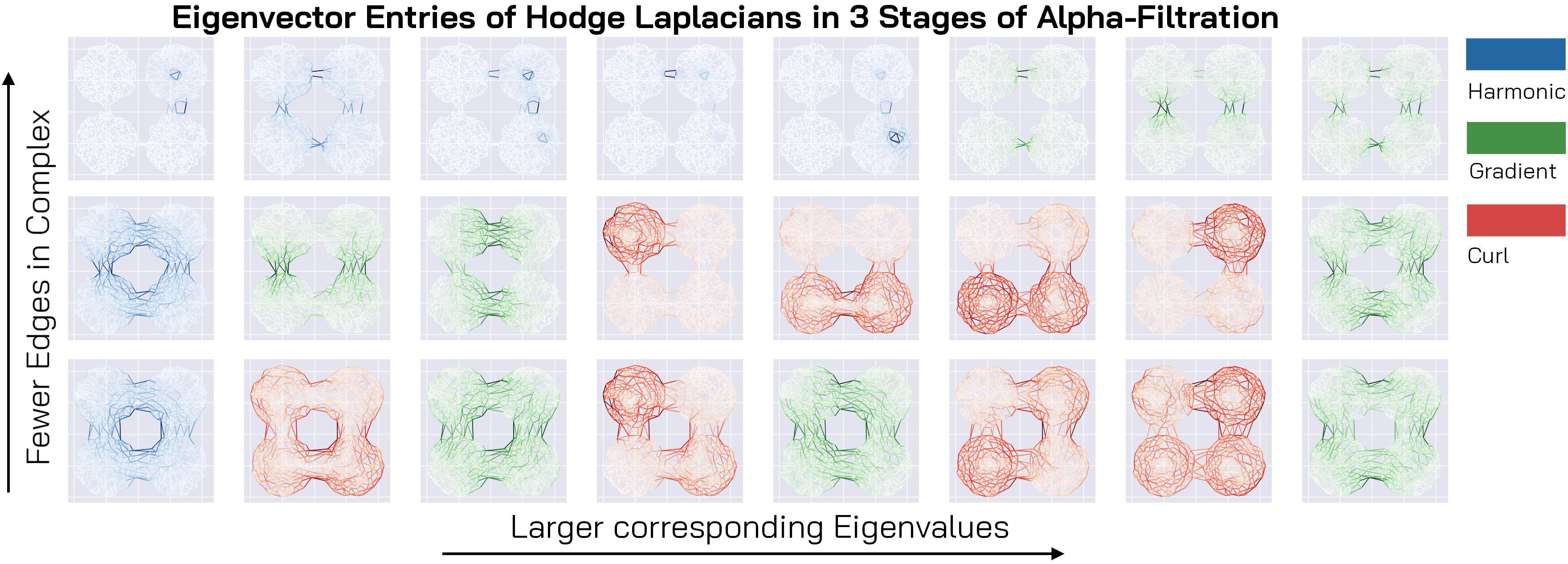}
			\vspace{-0.2in}
			\caption{\textbf{Smallest eigenvectors at three stages of the $\alpha$-filtration.} \emph{Blue:} Harmonic, \emph{Green:} Gradient, \emph{Red:} Curl eigenvectors.
				The intensity of the colour denotes the absolute value of the corresponding entry in the eigenvector, where we scale the maximal entry to 1 for visualisation purposes.
				The three stages of the filtration correspond to the rows and represent three different regimes:
				\textbf{1\textsuperscript{st} row:} There are many small \enquote{holes} in the data set. Accordingly, the first five eigenvectors are harmonic and strongly localised.
				In the subsequent eigenvalues, gradient components dominate, which are being localised on the bridges between the clusters.
				The curl eigenvalues are large in comparison.
				\textbf{2\textsuperscript{nd} row:} There is  a single harmonic eigenvector corresponding to the hole in the middle of the clusters, which is stronger on the bridges and the edges closer to the centre.
				The $\updownarrow$ and $\leftrightarrow$ gradient vectors are next and are localised at the bridges.
				We have a curl flow for each of the four clusters, and finally a diagonal gradient flow.
				\textbf{3\textsuperscript{rd} row:} There is a $1:1$ correspondence to the types of the previous vectors.
				However, the order of the eigenvectors has changed with the curl eigenvalues getting smaller.
				Both gradient and curl vectors are less localised.
			}
			\label{fig:eigenvectors}
		\end{center}
		\vspace{-0.3in}
	\end{figure*}
\section{Simplicial complexes and persistence}
\label{sec:background}
\subheading{Simplicial Complexes}
Simplicial complexes \cite{Hatcher:2002,Bredon:1993} are generalisations of graphs consisting of nodes ($0$-simplices), edges ($1$-simplices), filled-in triangles ($2$-simplices), tetrahedra ($3$-simplices), etc.
\begin{definition}[Simplicial complex]
	A simplicial complex $\SC$ consists of a set of vertices $X$ and a set of finite non-empty subsets of $X$, called simplices $S$, such that \textbf{(i)} $S$ is closed under taking non-empty subsets and
		\textbf{(ii)} for every $x\in X$, the singleton $\{x\}$ is contained in $S$.
	For simplicity, we often identify $\SC$ with its set of simplices and use $\SC_n$ to denote the subset of simplices with $n+1$ elements.
\end{definition}
The orientation of an simplex $\sigma \in \SC$ is represented by an ordering of its vertices.
We consider two orientations equivalent if they differ by an even permutation.
Hence, for $k>0$ every $k$-simplex admits two possible orientations.
We denote by $C_k\coloneqq \R^{|\SC_k|}$ the $k$-th simplicial signal space, i.e.\ the space of real signals supported on the $k$-simplices.
An orientation on the simplices allows for bookkeeping in form of \emph{boundary matrices}:
The boundary matrix $B_k(\SC)$ of an \tSC{} records the incidence relations between the $k$-simplices and the $(k-1)$-simplices with respect to the chosen orientations.
For $k=1$, this is just the ordinary edge-vertex incidence matrix of a graph.
We define the boundary matrix $B_0$ to be the empty matrix.

\subheading{Hodge Laplacians}
The Hodge Laplacians of an simplicial complex  $\SC$ are important linear maps with close connection to many properties of the underlying \textsmaller{SC}.
The $k$-th Hodge Laplacian $L_k$ is a square matrix indexed over the $k$-simplices of $\SC$.
\begin{definition}[Hodge Laplacian]
	For an \textsmaller{SC} $\SC$ with associated boundary matrices $\left(B_i\right)$, we denote by $L_k^\text{up}\coloneqq B_{k+1}B_{k+1}^\top$ the $k$-th Up-Laplacian, by $L_k^\text{down}\coloneqq B_k^\top B_k$ the $k$-th Down-Laplacian, and finally by their sum $L_k\coloneqq \DownL{k}+\UpL{k}$ the $k$-th \emph{Hodge Laplacian}. 
	We note that $L_0$ recovers the graph Laplacian.
\end{definition}
\begin{theorem}[Hodge Decomposition \cite{Lim:2020, Barbarossa:2020, Schaub:2021, Roddenberry:2021}]
	For an \tSC{} $\SC$ with boundary matrices $\left(B_i\right)$, Hodge Laplacians $\left(L_i\right)$, and simplicial signal spaces $\left(C_i\right)$, we have that
	\[
	\R^{|\SC_k|}\cong C_k =\lefteqn{\overbracket{\phantom{\underbrace{\Ima B_k^\top}_\text{gradient space} \oplus  \overbrace{\underbrace{\ker L_k}_\text{harmonic space}}}}^{\ker \UpL{k}}}
	\underbrace{\Ima B_k^\top}_\text{gradient space} \oplus  \overbracket{\underbrace{\ker L_k}_\text{harmonic space}\oplus \underbrace{\Ima B_{k+1}}_\text{curl space}}^{\phantom{\ker \DownL{k}}\ker \DownL{k}}.
	\]
\end{theorem}
For $k=1$, the \emph{harmonic} space characterises edge flows around the holes of the \tSC; the \emph{gradient} space characterises flows that can be generated by assigning \enquote{potentials} to every vertex and computing their difference along every edge; the curl space characterises flow \enquote{swirling} within the \tSC, and any curl flow can be created by assigning a (magnetic) potential to every 2-simplex in the \tSC{} and then computing the induced circulations around the boundaries of the 2-simplices. 
See \Cref{fig:eigenvectors} for a visualisation of different flow types. 
This structure has been exploited in many applications, with the harmonic space receiving the largest attention (\cite{Frantzen2021} for trajectory outlier prediction, \cite{Ebli2019} for simplicial clustering, and \cite{Grande:2023} for topological point cloud clustering).

\subheading{$\alpha$-Complexes and the Delaunay Triangulation}
Given a point cloud $X\subset \R^2$, a \emph{Delaunay} triangulation is a triangulation such that no point in $X$ is contained inside the circumcircle of any of the triangles.
We can extend this to an \tSC{} $\text{DT}(X)$ by taking $X$ as the vertices, the edges as the $1$-simplices, and the triangles as the $2$-simplices.
This construction naturally generalises to higher dimensions by dividing the convex hull of the points in $\R^n$ into $n$-simplices.
This construction provides a way to construct a sparse high-dimensional complex on top of data points in $\R^n$, on which we can then use ideas of topological signal processing.
However, we sometimes want to factor in a notion of \enquote{scale}, depending on which we might decide to connect only vertices which are \enquote{close}.
This leads to the notion of $\alpha$-complexes \cite{edelsbrunner2011alpha}, where $\alpha^2$ is essentially the threshold circumsphere radius $CR(\sigma)$ for a simplex $\sigma$ to be considered \enquote{small}:
\begin{definition}[$\alpha$-complexes]
	Given a point cloud $X\in \R^n$ with associated generalised Delaunay triangulation $\text{DT}(X)$ and a parameter $\alpha$.
	Then, the $\alpha$-complex is the \tSC{} $\SC_\alpha (X)$ with vertices $X$ and simplices 
$
	S_\alpha(X)=\left\{\sigma\in \text{DT}(X) :CR(\sigma)\le \alpha^2\right\}.
	$
\end{definition}
For $\alpha\le \alpha'$, we have the canonical inclusion $\SC_\alpha(X)\subset\SC_{\alpha'}(X)$.
Because of their origin in the Delaunay triangulation, the $\alpha$-complexes contain only very few edges in comparison to the popular Vietoris--Rips complexes, while retaining similar properties.
This results in
good computational performance in downstream tasks.

\subheading{Persistent Homology and the Persistent Laplacians}
In many applications, we have no a-priori knowledge of how to pick the scale-parameter $\alpha$ correctly.
Varying $\alpha$ in the $\alpha$-complexes of $X$ enables us to extract features of $X$ across all possible scales.
Different methods of defining, analysing, and subsuming these features have been proposed.
One of the most popular tools in this context is Persistent Homology, a tool from algebraic topology that tracks the evolution of the \enquote{shape} of a point cloud across the scales \cite{Edelsbrunner2008}, in terms of the number of connected components, holes, cavities, etc.
However, some interesting properties are not as black-and-white as the existence of holes:
For graphs, the second-smallest eigenvalue $\lambda_2$ of the graph Laplacian is a measure of how easily the graph can be divided into separate clusters~\cite{vonLuxburg:2007,chung1997spectral,spielman2012spectral}.
The extreme value of $0$ means that the graph is disconnected, while the maximum of $n$ is only achieved in complete graphs; all other graphs range somewhere in between.
Hence these \enquote{continuous} spectral quantities can contain more information than their binary counterparts from algebraic topology.
This has motivated the study of the persistent Laplacians \cite{wang2020persistent}, which are in a special case the Hodge Laplacians $L(\SC_\alpha (X))$ of the $\alpha$-complexes.
\section{Eigenvectors of the Hodge Laplacian}
\begin{figure}[tb!]
	\begin{center}
		\includegraphics[width=0.9\columnwidth]{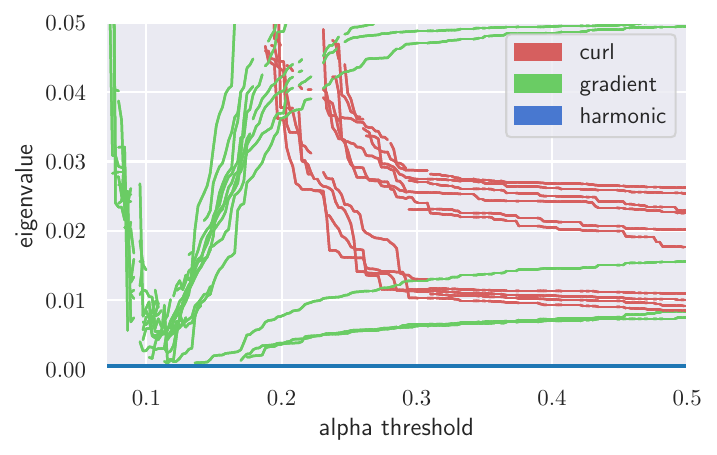}
		\includegraphics[width=0.9\columnwidth]{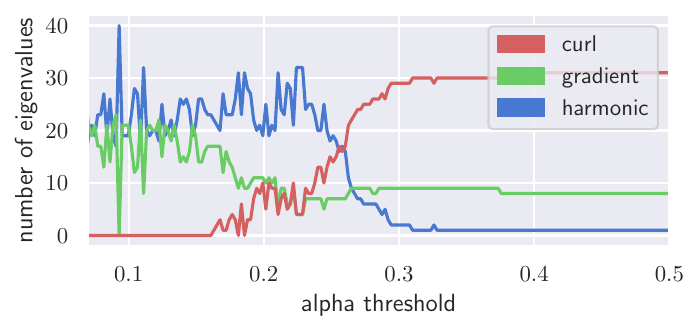}
		\vspace{-0.2in}
		\caption{\textbf{Evolution of the smallest eigenvalues of Hodge Laplacian across the steps of $\alpha$-filtration of the point cloud of \Cref{fig:eigenvectors}.}
			\emph{Top:} We witness two distinct phases in eigenvalue behaviour.
			In the \textbf{first regime} of low filtration values, the associated $\alpha$-complexes are sparsely connected and exhibit many holes (large number of number of harmonic eigenvalues $\lambda=0$).
			There is a V-shaped pattern of gradient eigenvalues.
            This correspond to the fact that the clusters first grow, which creates very ``thin'' bridges (small gradient eigenvalues) between growing clusters.
            Once the clusters start to grow together the connecting bridges become ``wider'' (more connecting edges), and the gradient eigenvalues increase again.
			In the \textbf{second regime} of medium to high filtration values, the associated $\alpha$-complexes grow denser and the holes are filled.
			Only one harmonic eigenvalue remains, whereas three structural gradient eigenvalues remain relatively small.
			However, the most distinct feature of the second regime is the cohort of decreasing curl eigenvectors/-values overtaking the gradient vectors of the initial V-shape.
			The smallest four of these eigenvalues are separated from the rest by a significant margin.
			They represent the curl flow generated around each of the four clusters.
			\emph{Bottom:} Number of eigenvalues in smallest $40$ being curl, gradient, and harmonic.
		}
		\label{fig:eigenvalues}
	\end{center}
	\vspace{-0.2in}
\end{figure}

\subheading{Harmonic, gradient, and curl eigenvectors}
\label{sec:manyeigenvectors}
In the previous section, we have introduced the concept of three different types of eigenvectors of Hodge Laplacians: harmonic, curl, and gradient eigenvectors.
In \Cref{fig:eigenvectors} we have depicted the eigenvectors corresponding to the smallest $8$ eigenvalues at three steps of an $\alpha$-filtration on a point cloud sampled from four disks arranged in a square.
It is immediately clear that a small eigenvalue can encode very different properties for the curl or the gradient case:
A small gradient eigenvalue arises from two clusters of the data connected by a small \enquote{bridge}.
A small curl eigenvalue arises from a densely connected cluster with large radius.
Finally, a harmonic eigenvalue just means that there exists a hole in the data set, without any information on its size.

\subheading{Tracking the eigenvectors across multiple filtration steps}
\Cref{fig:eigenvectors} reveals that there is a strong connection between eigenvectors in different steps of the filtration.
However, these eigenvectors equivalence classes frequently swap order.
Thus, simply treating the \emph{$n$-th smallest eigenvalue/vector pair} as a fixed instance does not respect the underlying structure of the filtration.
Accordingly, we propose to connect the eigenvectors based on similarity.
However, defining such a similarity is not straight-forward because the Hodge Laplacians of different filtration steps operate on different simplicial signal spaces.
Furthermore, being an eigenvector is invariant under multiplication by a scalar and hence two similar eigenvectors can differ by scalar multiplication by $-1$.

\begin{definition}[Persistent eigenvector similarity]
Let $X\in \R^n$ be a point cloud, $k\in \Z_{\ge 0}$ a degree, $0<\alpha<\alpha'$ two filtration values, and $L_{\alpha}$ and $L_{\alpha'}$ the $k$-th Hodge Laplacians of the $\alpha$-complexes $\SC_{\alpha}(X)$ and $\SC_{\alpha'}(X)$.
Furthermore, let $v$ and $v'$ be two eigenvectors of $L_{\alpha}$ and $L_{\alpha'}$ respectively.
Finally, let $\iota\colon \R^{|\SC_{\alpha,k}(X)|}\hookrightarrow \R^{|\SC_{\alpha',k}(X)|}$ be the inclusion of the simplicial signal spaces induced by the inclusion of $\alpha$-complexes $\SC_\alpha\hookrightarrow\SC_{\alpha'}$.
We then define the persistent eigenvector similarity (\textsmaller{PES}) of $v$ and $v'$ by
$
\text{PES}(v,v')\coloneqq |\iota(v)^\top v'|/\left(\lVert v\rVert_2\cdot \lVert v'\rVert_2\right).
$
\end{definition}
\begin{definition}[Persistent eigenvector matching (\textsmaller{PEM})]
	\label{def:PEM}
	Let $X\in \R^n$ be a point cloud, $k\in \Z_{\ge 0}$ a degree, $0<\alpha<\alpha'$ two filtration values, and $L_{\alpha}$ and $L_{\alpha'}$ the $k$-th Hodge Laplacians of the $\alpha$-complexes $\SC_{\alpha}(X)$ and $\SC_{\alpha'}(X)$.
	Furthermore, let $V$ and $V'$ be two sets of eigenvectors of $L_{\alpha}$ and $L_{\alpha'}$.
	We then define the \textsmaller{PEM}:
	\begin{align*}
	\text{PEM}_{X,\alpha,\alpha'}&(V,V')\coloneqq 	\smash{\left\{ (v,v')\in V\times V': \phantom{\min_{\hat{v}\in V} \text{PES} (\hat{v},v')= \text{PES}(v,v')=\min_{\hat{v}'\in V'}\text{PES} (v,\hat{v}')}\right.}\\
	&\left.\min_{\hat{v}\in V} \text{PES} (\hat{v},v')= \text{PES}(v,v')=\min_{\hat{v}'\in V'}\text{PES} (v,\hat{v}')\right\}
	\end{align*}
\end{definition}
For non-degenerate cases, i.e.\ when no two eigenvector pairs have the same \textsmaller{PES}, persistent eigenvector matching provides a one-to-one matching (possibly not matching some of the eigenvectors).
We can then use this matching to track the evolution of different eigenvalues across the stages of the filtration, see \Cref{fig:eigenvalues}.

\section{Applications}
\begin{algorithm}[tb]
	\caption{Hodge Spectral Clustering (\textsmaller{TSC})}
	\label{alg:TSC}
	\begin{algorithmic}
		\STATE {\bfseries Input:} Simplicial (\emph{/cellular/cubical}) complex $\SC$, dimension $n$,  \#eigenvectors $h$, \#clusters $c$, selection of  $C\in\{\text{gradient, harmonic, curl, total} \}$ eigenmodes.
		\STATE \textbf{1.} Compute smallest $h$ $C$-eigenvectors $v_1, \dots ,v_h$ of $L_n(\SC)$.
		\STATE \textbf{2.} Set $v(\sigma_n)\coloneqq \left(v_1\left[\sigma_n\right],\dots, v_h\left[\sigma_n\right]\right)$.
		\STATE \textbf{3.} Set $\hat{v}(\sigma_n)\coloneqq \text{sgn}\left(1^\top v(\sigma_n)\right)v(\sigma_n)$.
		\STATE \textbf{4.} Cluster $\SC_n$ using $k$-means on $\{\hat{v}(\sigma_n):\sigma_n\in \SC_n\}\subset \R^{|\SC_n|}$.
		\STATE \emph{Optional:} Cluster points in $\SC_0$ wrt.\ adjacent $n$-simplices.
		\STATE \textbf{Output:} Cluster-ids of $k$-simplices $\SC_n$ of step \textbf{4}.
	\end{algorithmic}
	\vskip -0.2in
\end{algorithm}

\begin{figure}[tb!]
	\begin{center}
		\begin{subfigure}{0.4\columnwidth}
			\includegraphics[width=\columnwidth]{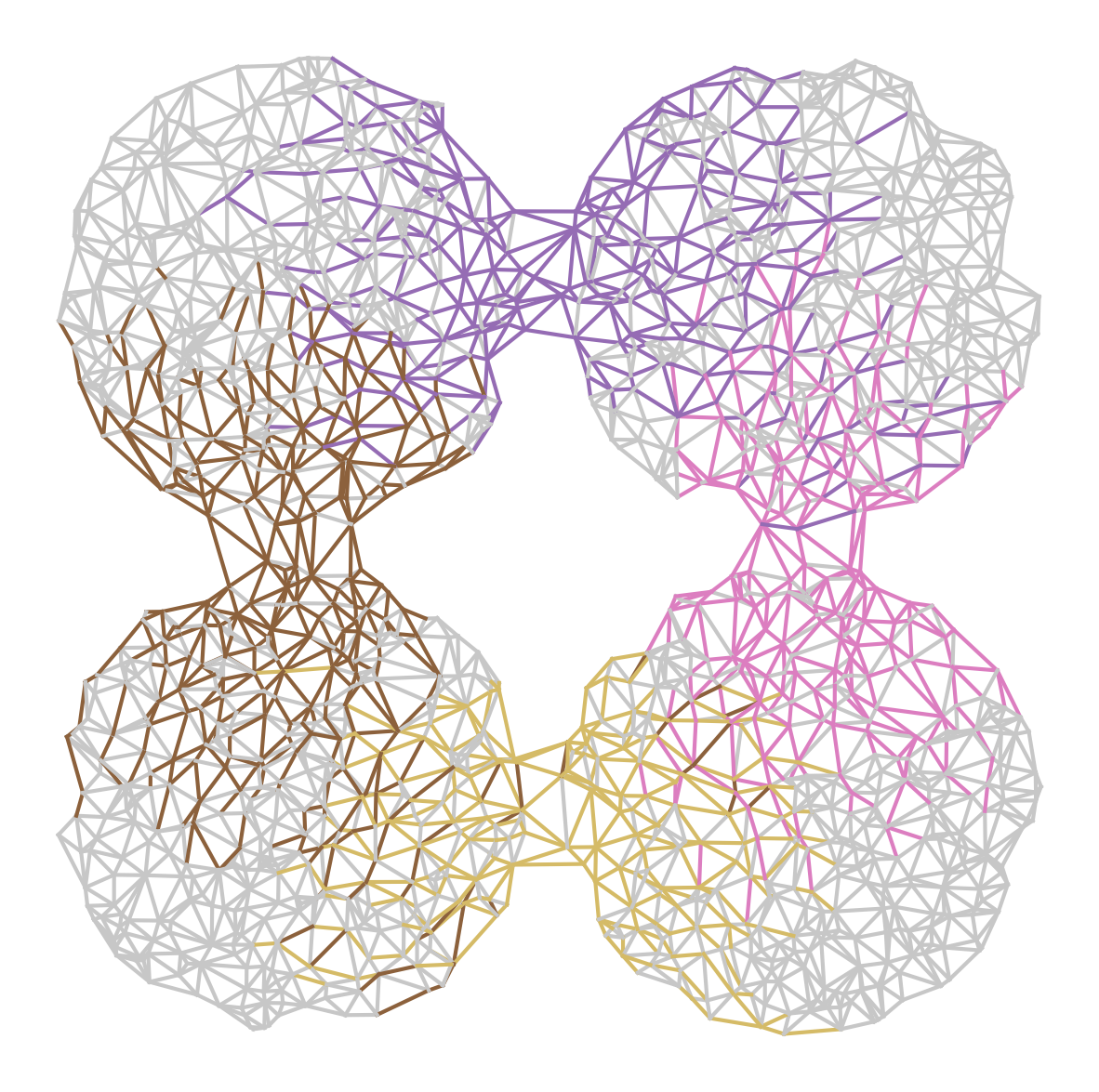}
			\subcaption{Gradient Clustering.}
		\end{subfigure}
		\begin{subfigure}{0.4\columnwidth}
			\includegraphics[width=\columnwidth]{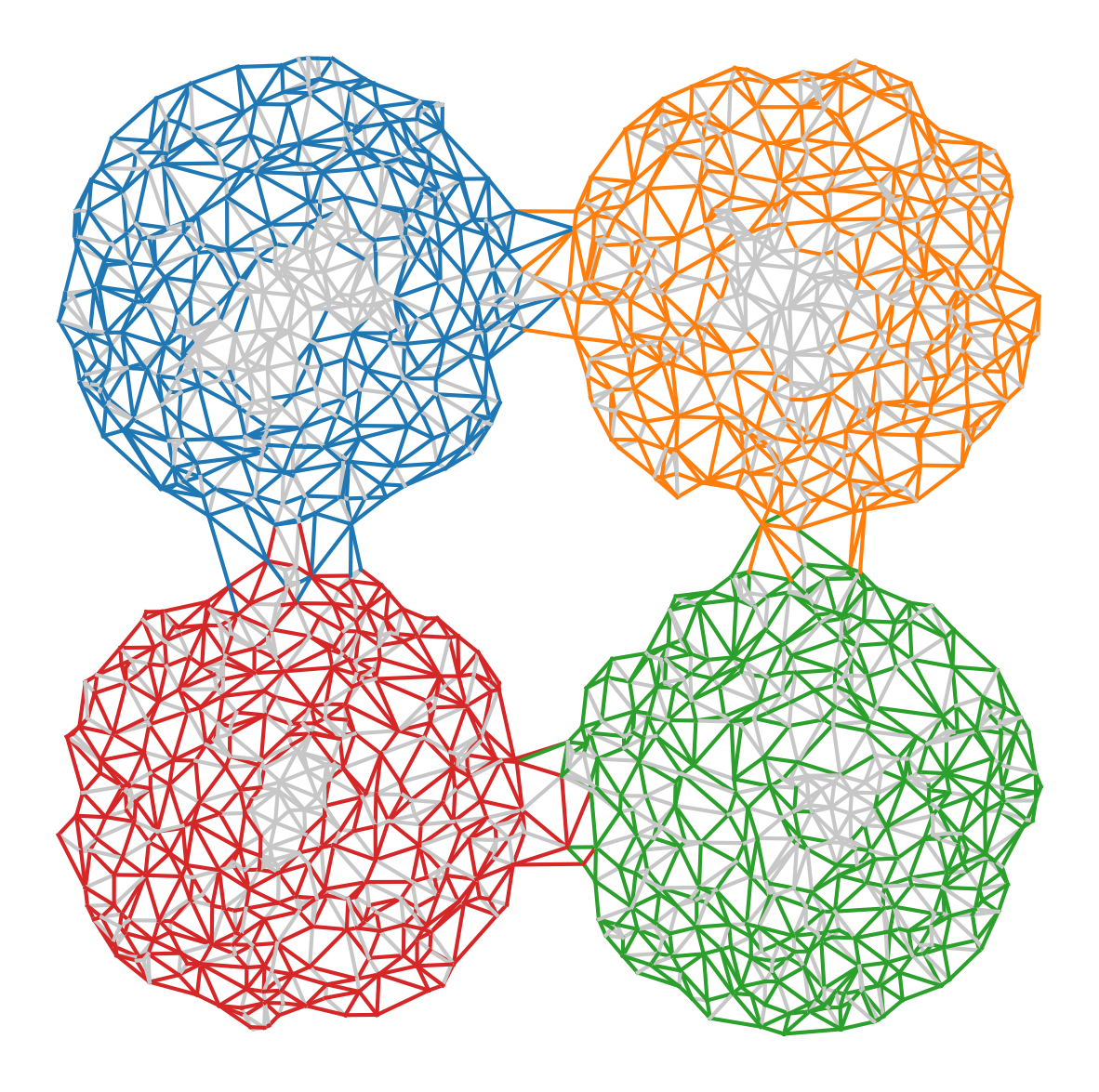}
			\subcaption{Curl Clustering.}
		\end{subfigure}
		\begin{subfigure}{0.4\columnwidth}
			\includegraphics[width=\columnwidth]{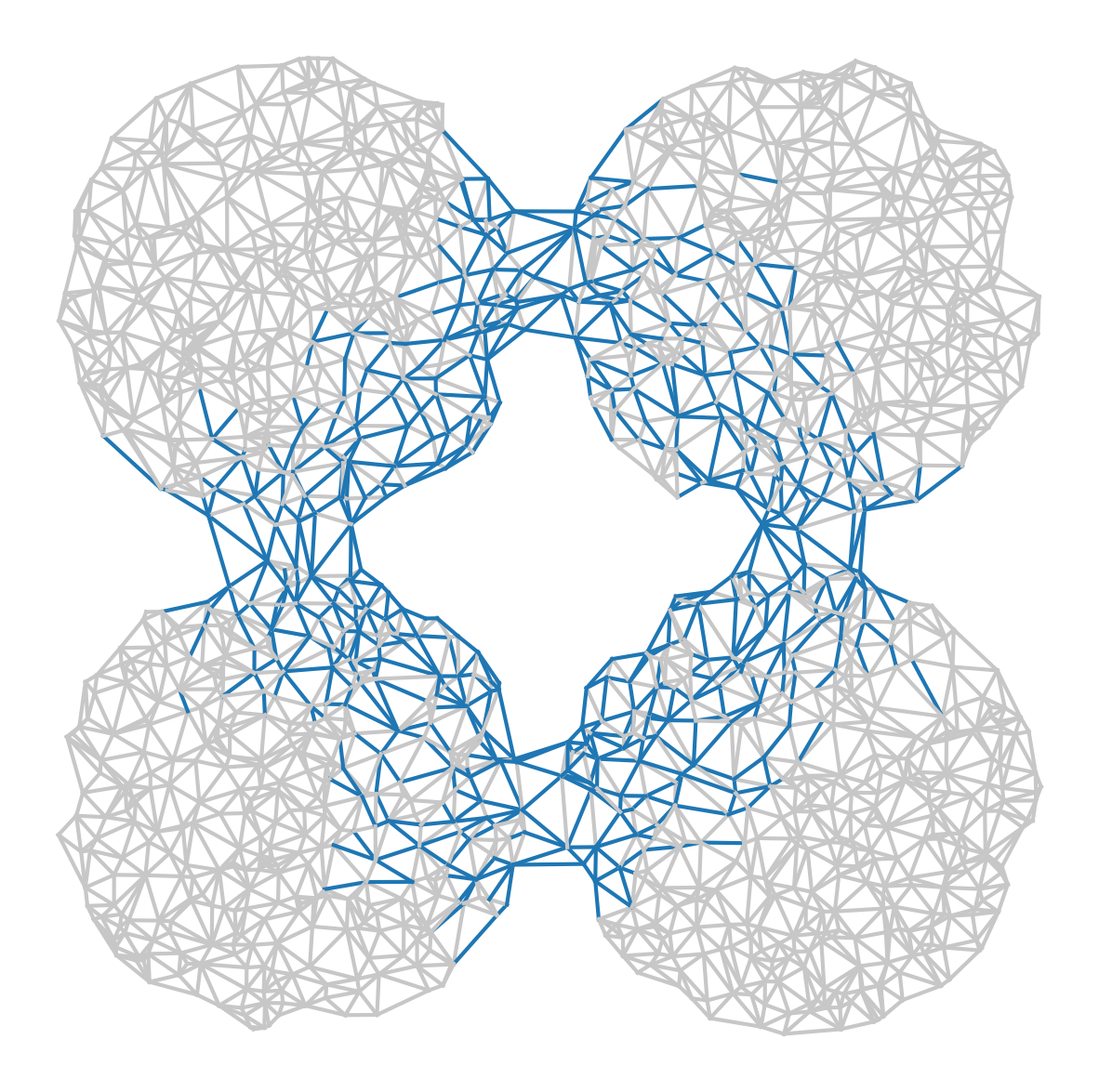}
			\subcaption{Harmonic Clustering.}
		\end{subfigure}
		\begin{subfigure}{0.4\columnwidth}
			\includegraphics[width=\columnwidth]{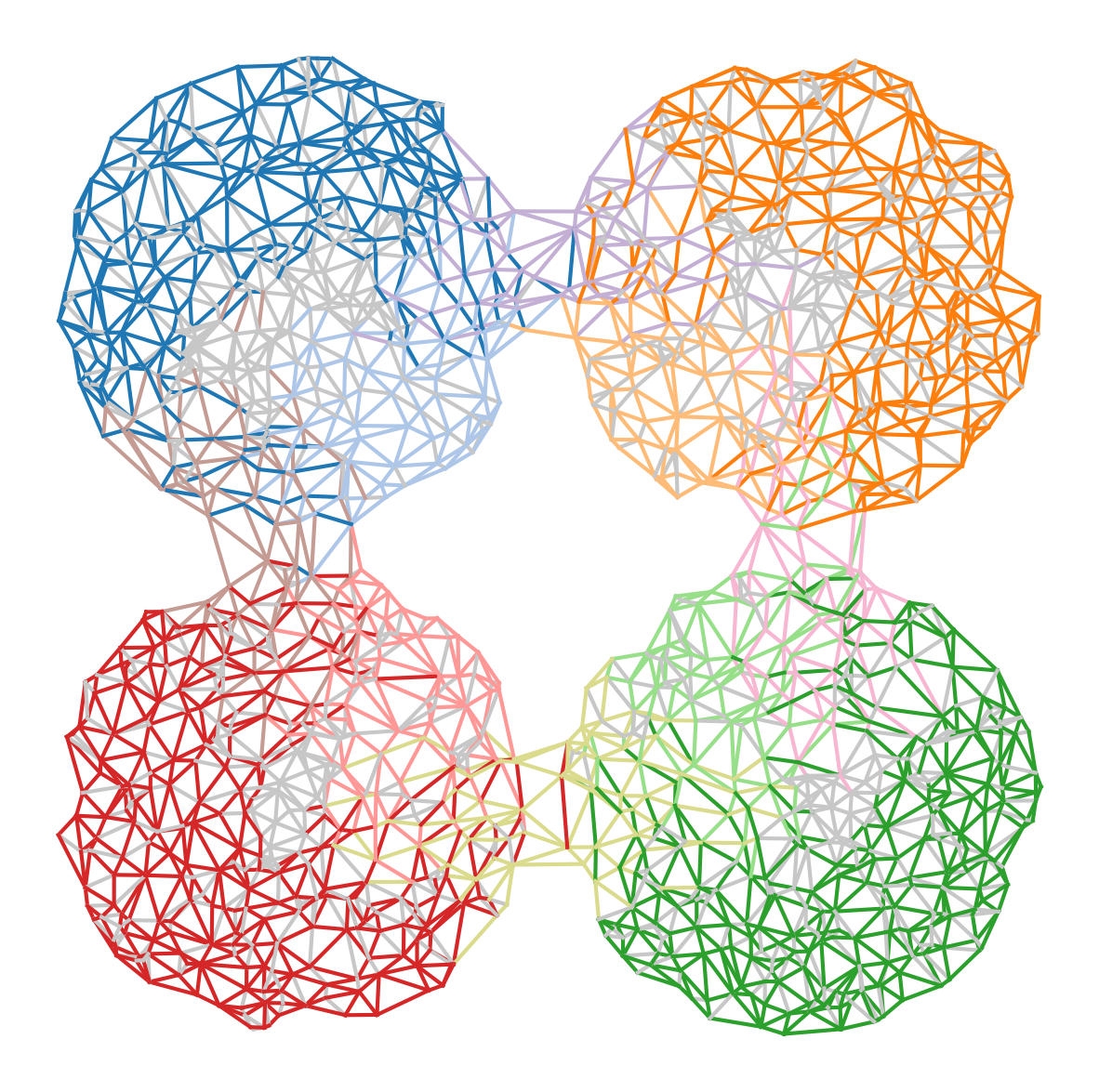}
			\subcaption{Total Clustering.}
		\end{subfigure}
		
		\caption{\textbf{Hodge Spectral Clustering.}
			In contrast to ordinary $0$-dimensional spectral clustering, higher-order spectral clustering allows for the subcategories of \emph{gradient}, \emph{curl} and \emph{harmonic} clustering, highlighting different features of the simplicial complex.
		}
		\label{fig:SpectralClustering}
	\end{center}
	\vspace{-0.3in}
\end{figure}

\label{sec:applications}
In the previous sections, we have argued why providing additional information on the often considered eigenvalues of Hodge Laplacians is crucial.
We now give a range of possible applications building on top of these considerations.

\subheading{Hodge spectral clustering}
We generalise the traditional notion of spectral clustering to higher-order simplicial complexes, see \Cref{alg:TSC}.
Given a simplicial complex $\SC$ with $n$-th Hodge Laplacian $L$ and eigenvectors $v_1,\dots,v_k$, we will cluster an $n$-simplex $\sigma_n\in \SC_n$ based on its associated eigenvector entries $v (\sigma_n)\coloneqq \left(v_1\left[\sigma_n\right],\dots, v_k\left[\sigma_n\right]\right)$.
However, in the construction of the Hodge Laplacian we have to fix an arbitrary orientation of the simplices, and changing this arbitrary orientation leads to sign changes in the eigenvector entries of $\sigma_n$.
Hence we need to cluster the set $\{v(\sigma_n):\sigma_n\in \SC_n \}/(v\sim -v)$ inside $\R^{|\SC_n|}$ where we identify $v$ and $-v$ for all $v\in \R^{|\SC_n|}$.
However, clustering algorithms work best on plain real vector spaces.
Thus we use $k$-means on
$\{\text{sgn}\left(1^\top v(\sigma_n)\right)v(\sigma_n):\sigma_n\in \SC_n\}$.
Note that we have made the embedding of $\sigma_n$ invariant under orientation change.
Building on our previous considerations, we can choose a subset of the curl, gradient, and harmonic eigenvectors to highlight different properties of our point cloud and our simplicial complex, see \Cref{fig:SpectralClustering}.
By aggregating the cluster-ids of the incident edges for each node, we can obtain an induced clustering of the points similar to~\cite{Grande:2023}.

\begin{figure}[tb!]
	\vspace{-0.1in}
	\begin{center}
		\begin{subfigure}{0.4\columnwidth}
			\includegraphics[width=\columnwidth, angle =180, origin=c]{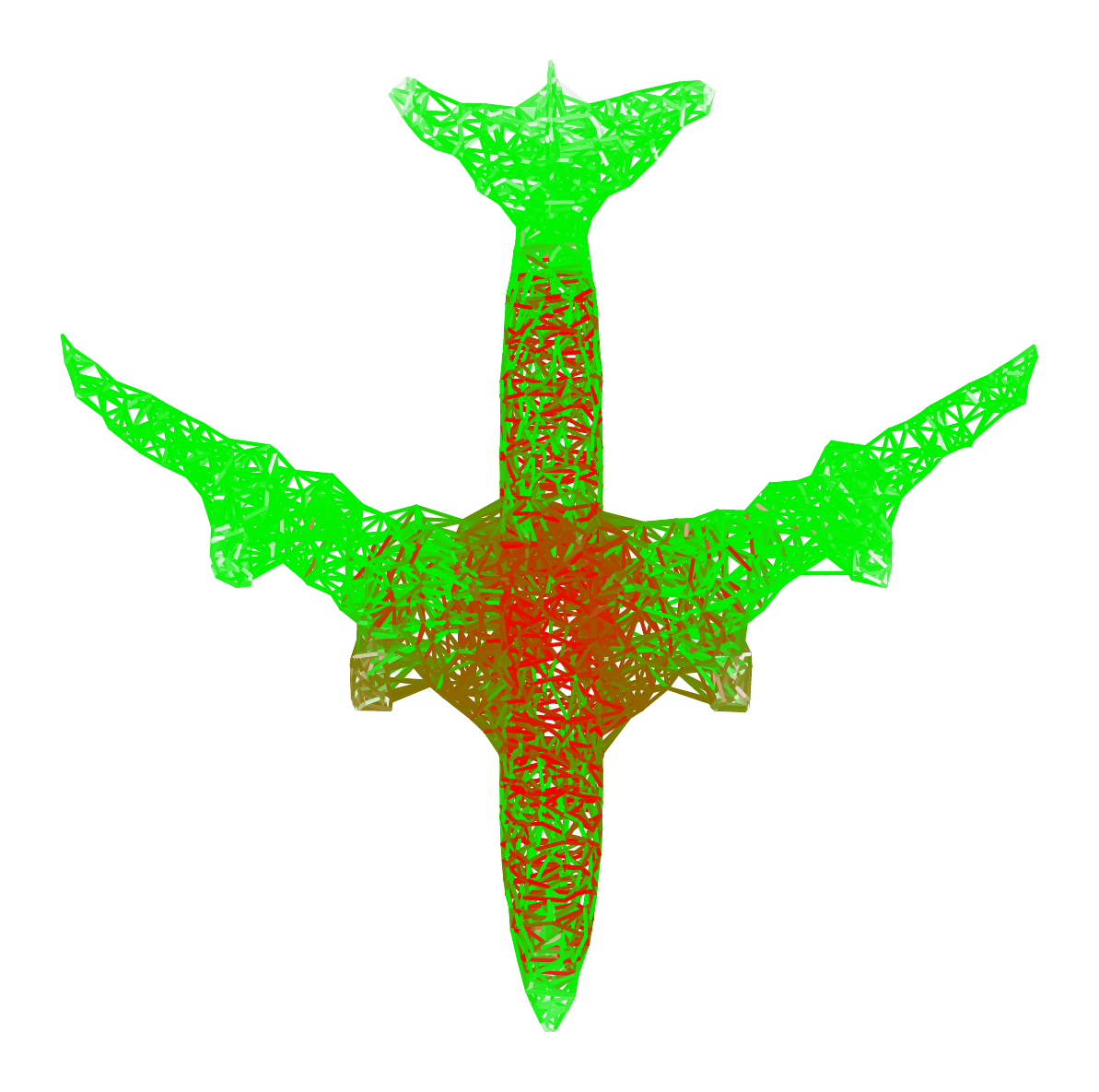}
		\end{subfigure}
		\begin{subfigure}{0.4\columnwidth}
			\includegraphics[width=\columnwidth]{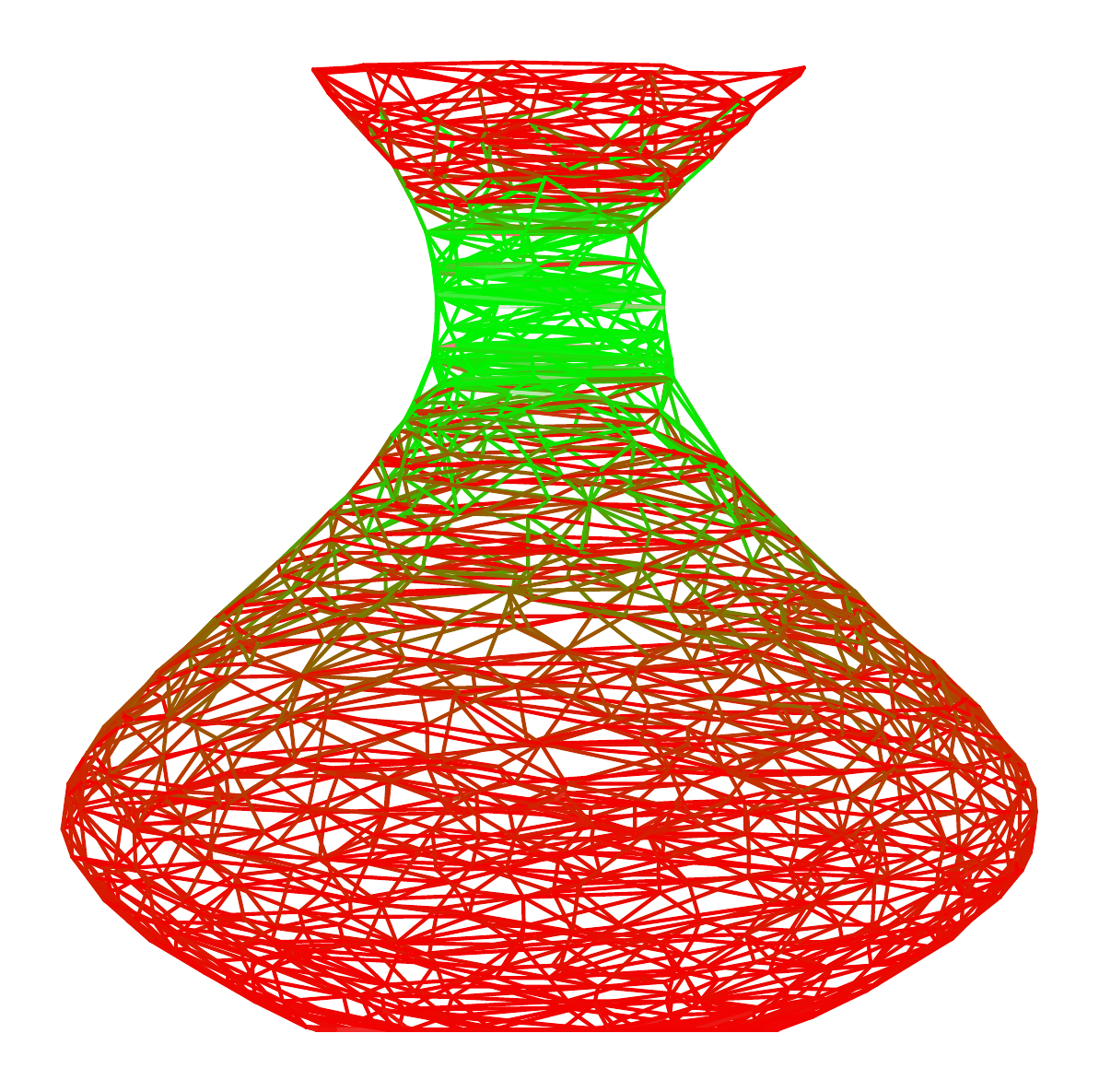}
		\end{subfigure}
		\begin{subfigure}{0.4\columnwidth}
			\includegraphics[width=\columnwidth]{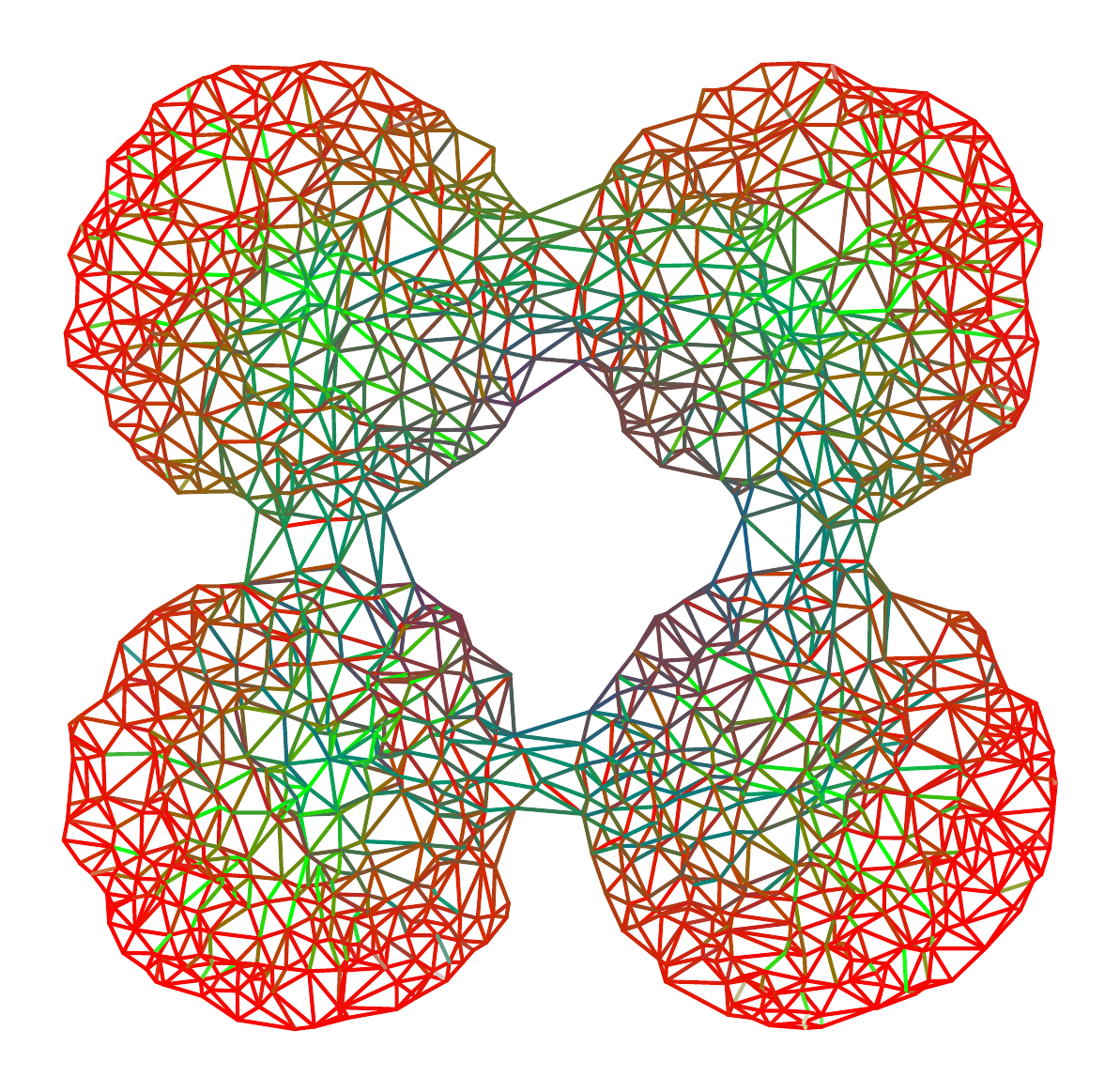}
		\end{subfigure}
		\begin{subfigure}{0.4\columnwidth}
			\includegraphics[width=\columnwidth]{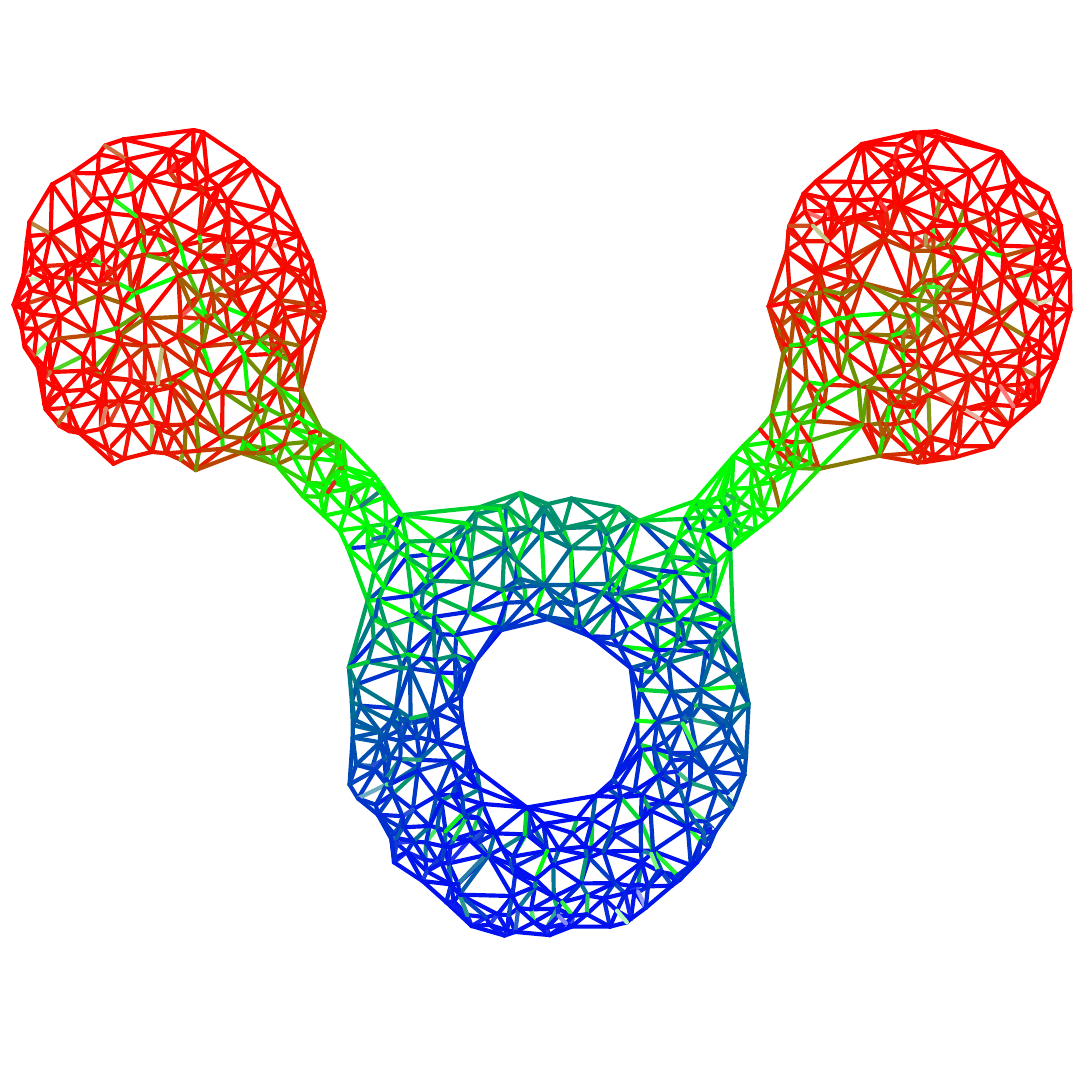}
		\end{subfigure}
		\vspace{-0.1in}
		\caption{\textbf{Edge roles using \textsmaller{HGC}-values.}
            Relative strength of \emph{blue:} harmonic, \emph{red:} curl, \emph{green:} gradient contributions. \emph{Top:} $\alpha$-complexes on two 3d point clouds sampled on real world objects of ModelNet40 data set \cite{wu20153d}. \emph{Bottom:} $\alpha$-complexes on two synthetic point clouds.
		}
		\label{fig:edgeroles}
	\end{center}
	\vspace{-0.2in}
\end{figure}
\subheading{Edge roles using \textsmaller{HGC}-values}
Edges ($n$-simplices) can have many different roles inside an \tSC:
Some edges connect nodes in the centre of a cluster, some are located in the outer rim of a cluster, or connect different clusters, and again others are located on the edge of \enquote{empty spaces} in the \tSC, or have no relationship to any of the clusters.
This role based clustering is orthogonal to \textsmaller{TSC}.
Edge role classification plays an important role in applications on brain networks in neuroscience~\cite{Faskowitz2022}, link analysis in social networks~\cite{demeo2012}, etc. 
We provide a metric to extract these relevant information:
\begin{definition}[\textsmaller{HGC}-values]
	Given a simplicial complex $\SC$, a dimension $n$ and a number $k$, let $V\coloneqq v_1,\dots, v_k$ be the eigenvectors associated to the $k$ smallest eigenvalues of the $n$-th Hodge Laplacian $L_n$ of $\SC$. Furthermore, we denote by $V_\text{harm}\coloneqq\{v_i\in V:L_nv_i = 0\}$, $V_\text{curl}\coloneqq\{v_i\in V:L^\text{up}_nv_i \not = 0\}$, $V_\text{grad}\coloneqq\{v_i\in V:L^\text{down}_nv_i \not = 0\}$ the harmonic, curl, and gradient eigenvectors.
	Let $e_\text{max}$ denote the entry with maximum absolute value across all $v\in V$.
	Setting the max over $\emptyset$  to be $0$, we can define the \textsmaller{HGC}-value of the simplex $\sigma_n$ $\text{HGC}(\sigma_n)$ to be
	\[
\left(\max_{v\in V_\text{harm}} (|v(\sigma_n)|) ,\max_{v\in V_\text{grad}} (|v(\sigma_n)|) ,\max_{v\in V_\text{curl}} (|v(\sigma_n)|) \right)/|e_\text{max}|.
\]
\end{definition}
Intuitively, the \textsmaller{HGC}-values measure the relevance of an edge for harmonic, gradient, and curl eigenvectors.
Combinations of these then correspond to the edge roles discussed above, see \Cref{fig:edgeroles} for an illustration.
Although normalising by the largest entry of any eigenvector worked best in practice, other normalisations
are possible as well.
For a different approach to edge roles based on harmonic, curl, and gradient vectors see \cite{schaub2014structure}.

\subheading{Discussion}
In this paper, we have argued why disentangling the roles and properties of different small eigenvalues of the Hodge Laplacian is beneficial in theory and applications.
Building on this insight, we have introduced a method to track individual eigenvalues through an $\alpha$-filtration, introduced a novel tunable spectral clustering algorithm on higher-order simplicial complexes, and introduced the notion of \textsmaller{HGC}-values to extract simplicial roles. Finally, we showed the expressiveness of Hodge spectral clustering and the \textsmaller{HGC}-values in experiments on synthetic and real-world data.
\section{Code}
The accompanying source code used to generate the experimental results of this paper can be found at \href{https://git.rwth-aachen.de/netsci/2024-disentangling-the-spectral-properties-of-the-hodge-laplacian}{this link}.
\bibliographystyle{IEEEbib}
\bibliography{refpc}
\end{document}